\documentstyle[12pt]{article}
\def\ba{\begin{array}}
\def\ea{\end{array}}    
\def\be{\begin{equation}}
\def\ee{\end{equation}}
\def\bem{\begin{em}}
\def\eem{\end{em}}
\def\ot{\otimes}

\def\a{\alpha}                           
\def\b{\beta}
\def\g{\gamma}
\def\D{\Delta}
\def\ep{\epsilon}
\def\ra{\longrightarrow}

\def\ca{{\cal A}}

\def\ca{{\cal A}}
\def\cb{{\cal B}}
\def\cc{{\cal C}}

\def\ch{{\cal H}}
\def\cm{{\cal M}}
\def\cn{{\cal N}}

\def\cp{{\cal P}}
\def\cu{{\cal U}}
\def\cv{{\cal V}}
\def\cw{{\cal W}}

\def\tg{\triangle}
\def\lb{\langle}
\def\rb{\rangle}
\newfont{\numb}{msbm10}
\def\com{\mbox{\numb C}}

\def\pg#1,#2,#3,{\langle #1 | #2 \rangle^{#3}}
\def\m#1,{m_{#1}}
\def\ic#1,{i_{#1}}
\def\jc#1,{j_{#1}}
\def\tt#1,{>\!\!\!\lhd_{#1}}
\def\tw{>\!\!\!\lhd}
\def\ps#1,#2,{\Psi_{#1,#2}}
\def\id#1,{id_{#1}}
\def\0n{_{(0)}}
\def\1n{_{(1)}}
\def\2n{_{(2)}}
\def\3{_{(3)}}
\def\ps#1,#2,{\Psi_{{#1}{,}{#2}}}
\def\pr{(.|.)}
\begin{document}
\title{On composite systems and quaqntum statistics}
\author{Wladslaw Marcinek\\
Institute of Theoretical Physics, University
of Wroc{\l}aw,\\ Pl. Maxa Borna 9, 50-204  Wroc{\l}aw,\\
Poland}
\maketitle
\begin{abstract}
An algebraic formalism for the study of a system of 
charged particles interacting with an external quantum field is
developed. The notion of monoidal categories with duality is
used for the description of composite systems and corresponding 
quantum statistics. The Fock space representation is also 
discussed.
\end{abstract}
\section{Introduction}
The concept of generalized quantum statistics and related
topics has been under intensive study over the past few years
\cite{owg,gre,moh,fi}. An algebraic formalism for particles with 
generalized statistics has been developed recently by the author 
in the series of papers \cite{WM3,WM4,WM6,WM8,wm7,mad,mco} and 
also \cite{qweyl,top,castat,sin,wmq,qsym}. In this attempt 
the creation and annihilation operators act on a quadratic algebra 
$\ca$. The creation operators act as the multiplication in $\ca$ 
and the annihilation ones act as a noncommutative contraction 
The algebra $\ca$ play the role of noncommutative Fock space. 
It is interesting that in this algebraic approach all commutation 
relations for particles equipped with arbitrary statistics can 
be described as a representation of the so--called quantum 
Weyl algebra $\cw$ (or Wick algebra) \cite{jswe,RM,m10,ral}. 
The applications for particles in singular magnetic field 
has been given in \cite{mco,top,sin}. 
Similar approach has been also considered by others authors, 
see \cite{sci,twy,mep} and \cite{mphi,mira} for example.
An interesting approach to quantum statistics has been 
also given in \cite{fios,melme}.

In this paper an algebraic model of composite systems and 
quantum statistics is considered. The starting point for 
the study of our model is a system of charged particles 
interacting with certain external quantum field. The proper 
physical nature of the system is not essential for our 
considerations.  We assume that the system of charged particles 
before interaction are described by a monoidal (tensor) category 
$\cc$. Different species of quanta of the external field are 
described by a finite Hopf algebra $H$. We also assume that every 
charged particle is equipped with ability to absorb and emit quanta
of the external field. A system which contains a charge and certain 
number of quanta as a result of interaction with the external field 
is said to be a composite system of quasiparticles. A quasiparticle 
is a charged particle dressed with a single quantum of the external 
field. Note that in our approach a charged particle equipped with two 
quanta of the external field is considered as a system of two 
quasiparticles. Two quasiparticles are said to be identical if they 
are dressed with quanta of the same species. In the opposite case 
when the particle is equipped with two different species of quanta 
then we have different quasiparticles. It is interesting that 
quasiparticles have also their own statistics. Next we assume 
that our system after interactions is described by an another 
monoidal category $\cc^{ext}$, an extension of $\cc$. 
\section{Preliminaries}
In this paper by $k$ we denote in general an arbitrary field, but 
for the physical applications we restrict our attention to the the 
field of complex numbers $k \equiv \com$. First let us briefly 
recall the concept of monoidal categories. A {\it monoidal category} 
$\cc \equiv \cc(\ot, k)$ is shortly speaking a category equipped 
with a monoidal associative operation (a generalized
tensor product) $\ot :{\cal C}\times{\cal C}\ra {\cal C}$ which 
has a two-sided identity object $k$, see 
\cite{Maj,bm,qm,ML,jst,jet,man} for details. 
A {\it category with duality}
$\cc \equiv \cc(\ot, k, \ast, g, \Psi)$ is a monoidal 
category $\cc(\ot, k)$ equipped in addition with a 
$\ast$--operation, a pairing $g$, and a cross symmetry $\Psi$. 
The existence of $\ast$-operation means that for every object $\cu$ 
of the category $\cc$ there is an object $\cu^{\ast}$ of $\cc$ called
the {\it dual} of $\cu$. The dual $\cu^{\ast}$ of $\cu$ is assumed to 
be 
unique. Moreover we assume that the second dual $\cu^{\ast\ast}$ 
of $\cu$ is identical with $\cu$ and 
$(\cu \ot \cv)^{\ast} = \cv^{\ast} \ot \cu^{\ast}$ for every 
two objects $\cu$ and $\cv$ of $\cc$. If $f : \cu \ra \cv$ 
is a morphisms in $\cc$, 
then there is also a corresponding {\it dual morphisms} 
$f^{\ast} : \cv^{\ast} \ra \cu^{\ast}$. 
A (right) pairing in the category $\cc$ is a set of maps
\be
\ba{c}
g \equiv \{g_{\cu} : \cu^{\ast} \ot \cu \ra k \}
\label{parc}
\ea
\ee
satisfying some compatibility axioms \cite{castat}.
The cross symmetry is a set of natural isomorphisms 
\be
\Psi \equiv \{\ps \cu^{\ast}, \cv, : 
\cu^{\ast} \ot \cv \ra \cv \ot \cu^{\ast} \}.
\ee
such that
\be
\ba{l}
\Psi_{\cu^{\ast} \ot \cv^{\ast} ,\cw} 
= (\Psi_{\cu^{\ast},\cw} \ot \id \cv,) 
\circ (\id \cu, \ot \Psi_{\cv^{\ast} ,\cw}),\\
\Psi_{\cu^{\ast} ,\cv \ot \cw} 
= (\id \cv, \ot \Psi_{\cu^{\ast} ,\cw}) 
\circ (\Psi_{\cu^{\ast} ,\cv} \ot \id \cw,),
\label{heb}
\ea
\ee
for every objects $\cu, \cv, \cw$ in $\cc$ \cite{castat}. 
Note that the generalized cross symmetry is not a braid symmetry
in general. For a braid symmetry in the category with duality
$\cc$ we need the following additional two sets of operators
\be
\ba{l}
\Psi \equiv \{\ps \cu, \cv, : \cu \ot \cv \ra \cv \ot \cu\}\\
\Psi \equiv \{\ps \cu^{\ast}, \cv^{\ast}, : 
\cu^{\ast} \ot \cv^{\ast} \ra \cv^{\ast} \ot \cu^{\ast} \}
\ea
\ee
where $\cu ,\cv$ and $\cu^{\ast} , \cv^{\ast}$ are objects
of the category. We need also some new commutative diagrams for 
all these transformations and pairings \cite{Maj,jst,jet,ML}. 
The monoidal category with unique duality and braid symmetry 
is said to be rigid \cite{Maj,bm,qm}. It is denoted by 
$\cc \equiv \cc(\ot, k, \Psi)$. If in addition $\Psi_{U,V}^2 = id$, 
for every objects $U, V \in \cc$, then the set of all such
operation is said to be a symmetry or tensor symmetry and it
is denoted by $S$ instead of $\Psi$. The corresponding category $\cc$
is called a symmetric monoidal or tensor category \cite{man,Lub,GRR}.
An algebra $\ca$ equipped with a multiplication 
$m : \ca \ot \ca \ra \ca$ and unit 
$u : k \ra \ca$ which is at the same time an object of the 
category $\cc$ is said to be an algebra in the category. 
For the construction of the category $\cp$ corresponding to 
the system of charged particles before interactions we assume 
that there is:

$\bullet$ a linear space $E$ called a space of elementary states,

$\bullet$ the complex conjugated space $E^{\ast}$ for the

description of conjugated states,

$\bullet$ a composition operation $\ot$ which tell us how to

built composite states from elementary ones,

$\bullet$ a pairing $\pr_E : E^{\ast} \ot E \ra \com$ for the

construction of scalar product,

$\bullet$ a linear operator $C : E^{\ast}\ot E\ra E\ot E^{\ast}$ 

called an $E$--cross operator which describes the statistics of

the initial particles.\\
The category $\cp$ contains: the basic field $\com$, the space $E$ 
and the conjugate space $E^{\ast}$, all their tensor products and
direct sums. The pairing $\pr_E$ and the $E$--cross operator $S$ 
are initial objects for the construction of the pairing $g$ and 
the cross symmetry $\Psi^C$ in the category $\cp$. We also assume
that there is an algebra $\Lambda$ in the category $\cp$ which 
describes quantum states of charged particles. If the category
$\cp$ is symmetric monoidal, then we assume in addition that the
algebra is $S$-symmetric, i.e. $m = m \circ S$, where $m$ is the
multiplication in $\Lambda$ and $S \equiv \Psi^C$.
For the description of interactions of our system of charged 
particles with an external field we assume that there is a 
finite Hopf algebra $H$. The Hopf algebra $H$ describes 
different species of quanta of an external field. 
The system after interactions is described by an extended
category $\cp^{ext}$. This category describes  quantum states 
of composite systems of quasiparticles. 
Let $M$ be an arbitrary object of the category $\cp^{ext}$.
We assume that there is a coaction $\rho_E$ of $H$ on $M$, 
i.e. a linear mapping $\rho : H \ra M \ot H$, which define a 
(right-) $H$-comodule structure on $M$. The 
coaction $\rho$ of $H$ on $M$ describes the ability of
particles to the process of absorption of quanta of the external 
field. If we want to describe at the same time the process of
emission we need also a module action of $H$ on $M$.
Obviously these two structures, i. e. comodule and module
structure on $M$ must be compatible. In this way we obtain
that $M$ is the so-called right $H$-Hopf module. We also
obtain that our category $\cp^{ext}$ should be in fact a 
category of right $H$-Hopf modules. 
For the construction of the algebraic Fock space representation
we also need the notion of duality in our categories.
For this goal we want the notion of left $H^{\ast}$-modules,
comodules and Hopf modules, where $H^{\ast}$ is the dual Hopf 
algebra with respect to a pairing $g_H : H^{\ast}\ot H \ra k$. 
We also need a linear operator
$T : H^{\ast} \ot H \ra H \ot H^{\ast}$ called the $H$--cross
operator. This operator describes the statistics of quanta
of the external field. One can assume in addition that $H$ 
is selfdual. Our final category is the category 
$\cp^{ext} := \ _{H^{\ast}} ^{H^{\ast}}\cp _H^H$ of right 
$H$--Hopf modules and left $H^{\ast}$-Hopf modules. 
The statistics of composite systems is described by the cross 
symmetry in the category $\cp^{ext}$. Quantum states of the 
system after interactions are described by an algebra $\ca$
in the category $\cp^{ext}$. The algebra $\ca$ is an extension
of the algebra $\Lambda$ in $\cp$.
\section{Algebras, coalgebras and modules}
An unital and associative algebra is a $k$-linear
space $\ca$ equipped with two $k$-linear maps: a multiplication
$m : \ca \ot \ca \ra \ca$ and an unit $u : k \ra \ca$ such that
the following diagrams are commutative
\be
\ba{rccccl}
\ca \ot \ca \ot \ca&\ot&\ca&\stackrel{m\ot id}\ra&\ca \ot \ca&\\
&&&&&\\
\scriptstyle{id \ot m}&\downarrow&&&\downarrow&\scriptstyle{m}\\
&&&&&\\
\ca&\ot&\ca&\stackrel{m}\ra&\ca&\\
\label{alg1}
\ea
\ee
\be
\ba{rcccl}
&&\ca \ot \ca&&\\
\scriptstyle{u\ot id}&&&&\scriptstyle{id \ot u}\\
&\nearrow&|&\searrow&\\
&&|&&\\
k\ot\ca&&|&&\ca\ot k\\
&\searrow&\downarrow&\swarrow&\\
&&&&\\
&&\ca&&
\label{alg2}
\ea
\ee
A coalgebra is a $k$-linear space $\cc$ equipped with 
two $k$-linear maps: a comultiplication $\Delta : \ca \ot \ca \ra 
\ca$ 
and a counit $\eta : \ca \ra k$ such that
the following diagrams are commutative
\be
\ba{rcccl}
&&\Delta&&\\
&\cc&\ra&\cc \ot \cc&\\
&&&&\\
\Delta&\downarrow&&\downarrow&\Delta \ot id\\
&&&&\\
&\cc \ot \cc&\ra&\cc \ot \cc \ot \cc \ot \cc&\\
&&id \ot \Delta&&
\label{coalg1}
\ea
\ee
\be
\ba{rcccl}
&&\cc&&\\
&&&&\\
&\swarrow&|&\searrow&\\
&&|&&\\
k\ot\cc&&|&&\cc\ot k\\
&&|&&\\
\eta\ot id&\nwarrow&\downarrow&\nearrow&id\ot\eta\\
&&&&\\
&&\cc \ot \cc&&
\label{coalg2}
\ea
\ee
A bialgebra is a $k$-linear space $\cb$ equipped with both
algebra $\cb (m, u)$ and coalgebra $\cb (\D, \eta)$ structure
such that $\D$ and $\eta$ are algebra morphisms or equivalently,
$m$ and $u$ are coalgebra morphisms.
A Hopf algebra is a $k$-linear space $H$ equipped with a bialgebra
$H(m, u, \D , \eta)$ structure and an element $S$ called an atipode,
i. e. a linear map $S : H \ra H$ such the diagram
\be
\ba{lcccccr}
&&H\ot H&\ra&H\ot H&&\\
&&&&&\\
\tg&\nearrow&&&&\searrow&m\\
&&&&&&\\
H&&\ra&k&\ra&&H\\
&&&&&&\\
\tg&\searrow&&&&\nearrow&m\\
&&&&&&\\
&&H\ot H&\ra&H\ot H&&
\label{anty}
\ea
\ee
is commutative. 

Let $\cc$ be a $k$-coalgebra. A (right) $\cc$-comodule is a $k$-linear
space $\cn$ equipped a $k$-linear map $\Phi : \cn \ra \cn \ot \cc$
such that the following diagrams are commutative
\be
\ba{rccllrcccl}
&\cn&\stackrel{\Phi}\ra&\cn \ot \cc& \quad
&\cn&\stackrel{\Phi}\ra&\quad \cn \ot \cc&\\
&&&&\quad &&&&\\
\Phi&\downarrow&&\downarrow id\ot\tg& & \quad
&\searrow&\downarrow&id \ot \eta\\
&&&& \quad&&&&\\
&\cn \ot \cc&\stackrel{\Phi \ot id}\ra&\cn \ot \cc \ot \cc&
\quad &&&\cn \ot k&
\label{comod1}
\ea
\ee
Let $H$ be a Hopf algebra over o field $k$. We use the 
following notation for the coproduct in $H$: if $h \in H$, 
then $\tg (h) := \Sigma h\1n \ot h\2n \in H \ot H$. 
A {\it coquasitriangular Hopf algebra} (a CQTHA) 
is a Hopf algebra $H$ equipped with a bilinear form 
$\lb -,-\rb : H \ot H \ra k$ such that
\be
\ba{l}
\Sigma \lb h\1n, k_1 \rb k\2n h\2n = 
\Sigma h\1n k\1n \lb h\2n, k\2n \rb ,\\
\lb h, kl \rb = \Sigma \lb h\1n, k \rb \lb h\2n, l \rb,\\
\lb hk, l\rb = \Sigma \lb h, l\2n \rb \lb k, l\1n \rb
\ea
\ee
for every $h, k, l \in H$. If such bilinear form $b$ exists
for a given Hopf algebra $H$, then we say that there is a 
{\it coquasitriangular structure} on $H$.

Let $H$ be a CQTHA with coquasitriangular structure 
$\lb -,-\rb$. The family of all $H$-comodules forms 
a category $\cc = \cm^H$. The category $\cc$ is braided 
monoidal. The braid symmetry 
$\Psi \equiv \{\ps U, V, : U\ot V\ra V\ot U; U, V\in Ob\cc\}$
in $\cc$ is defined by the equation
\be
\ba{c}
\ps U, V, (u \ot v) = 
\Sigma\lb v\1n , u\1n \rb \ v\0n \ot u\0n ,
\label{coin}
\ea
\ee
where $\rho (u) = \Sigma u\0n \ot u\1n \in U \ot H$, and
$\rho (v) = \Sigma v\0n \ot v\1n \in V \ot H$ for every 
$u \in U , v \in V$.

Let $G$ be an arbitrary group, then the group algebra $H := kG$
is a Hopf algebra for which the comultiplication, the counit,
and the antipode are given by the formulae
$$
\ba{cccc}
\tg (g) := g \ot g,&\eta(g) := 1,&S(g) 
:= g^{-1}&\mbox{for} \ g \in G.
\ea
$$
respectively. If $H \equiv kG$, where $G$ is an Abelian 
group, $k \equiv \com$ is the field of complex numbers, 
then the coquasitriangular structure on $H$ is given
as a bicharacter on $G$ \cite{mon}. Note that for
Abelian groups we use the additive notation. A mapping 
$\ep : G \times G \ra \com \setminus \{0\}$ is said to 
be a {\it bicharacter} on $G$ if and only if we have 
the following relations
\be
\ep (\a, \b + \g) = \ep (\a, \b) \ep (\a , \g), \quad
\ep (\a + \b , \g) = \ep (\a , \g) \ep (\b , \g)
\ee
for $\a, \b, \g \in G$. If in addition 
$\ep (\a, \b) \ep (\b, \a) = 1$, for $\a, \b \in G$, 
then $\ep$ is said to be a {\it normalized bicharacter}. 
Note that this mapping are also said to be a commutation 
factor on $G$ and it has been studied previously by 
Scheunert \cite{sch} in the context of color Lie algebras.
We restrict here our attention to normalized bicharacters 
only. It is interesting that the coaction of $H$ on certain 
space $E$, where $H \equiv kG$ is equivalent to the 
$G$-gradation of $E$, see \cite{mon}.

If $\rho : M \ra M \ot H$ is a right $H$-comodule map, the set
$M^{coH}$ of $H$-coinvarints is defined by the formula
\be
M^{coH} := \{m \in M: \rho(m) = m \ot 1\}.
\ee
This means that $m^{coH}$ is trivial right $H$-comodule.
If $E$ is a trivial right $H$-comodule, then $E \ot H$
is a nontrivial right $H$-comodule. The comodule map
$\rho : E\ot H \ra E \ot H \ot H$ is given by the relation
\be
\rho(x \ot h) := \Sigma \ x \ot h\1n \ot h\2n,
\ee
where $x \in E, h \in H$ and 
$\tg (h) = \Sigma \  h\1n \ot h\2n$. The inclusion of the 
notion of right $H$-modules leads to the concept of Hopf 
modules \cite{mon}. A right $H$-Hopf module is a $k$-linear 
space $M$ such that (i) there is a right $H$-module
action $\lhd : M \ot H \ra M$, (ii) there is a right
$H$-comodule map $\rho : M \ra M \ot H$, (iii) $\rho$ is
a right $H$-module map, this means that we have the relation
\be
(m \lhd h)\0n \ot (m \lhd h)\1n = \Sigma \
m\0n \lhd h\0n \ot m\1n \lhd h\2n,
\ee
where $m \in M, h \in H, \rho (m) = \Sigma \ m\0n \ot m\1n$,
and $\tg (h) = \Sigma \  h\1n \ot h\2n$. It is very interesting 
that every right $H$-Hopf module $M$ is a tensor product 
$M = \cu \ot H$, where $\cu \equiv M^{coH}$ is also a trivial
right $H$-module, i. e. $m \lhd h = \eta (h) m$ for every 
$m \in M$. If $H = kG$, then every right $H$-Hopf module $M$
is a $G$-graded space $M = \oplus_{g\in G}M_g$, such that
$M_g := \cu \ot g$, where $g \in G$ and $\cu$ is an arbitrary
linear space with trivial action of $G$. In this case we have
\be
m_g \lhd h = m_{gh}, \;\;\; \rho (m_g) = m_g \ot g,
\ee
where $m_g \in M_g$, $m_g := m \ot g$, $m \in \cu$, $g, h \in G$.

Let $\cp$ be a monoidal category such that $\cu \equiv \cu^{coH}$ 
for every object $\cu$ in $\cp$, then by $\cp^H$ we denote the
monoidal category of right $H$-comodules of the form 
$M := \cu \ot H$, where $\cu$ is an object in $\cp$, and $H$ 
is a Hopf algebra. The category $\cp^H$ is said to be a right 
$H$-comodule extension of $\cp$. If in addition every object
$\cu$ of $\cp$ is a trivial right $H$-module, then we obtain
the category of right $H$-Hopf modules $\cp^H_H$.

Let $H^{\ast}$ be a Hopf algebra, the dual of $H$. Then we
introduce the notion of left $H^{\ast}$-modules, comodules,
Hopf modules and corresponding categories in a similar way.
In this way we obtain the category 
$\cp^{ext} := \ _{H^{\ast}} ^{H^{\ast}}\!\cp _H^H$ of right 
$H$-Hopf modules and left $H^{\ast}$-Hopf modules. The
pairing $g^{ext}$ in $\cp^{ext}$ is given by the formula
\be
\ba{c}
g^{ext}_M ((h^{\ast} \ot m^{\ast}) \ot (m \ot h)) := 
g_{\cu}(m^{\ast}\ot m) g_H (h^{\ast}\ot h),
\label{epa}
\ea
\ee
where $M := \cu \ot H$, $\cu$ is an objects in $\cp$, 
$M^{\ast} := H^{\ast} \ot \cu^{\ast}$, $m \in M$, 
$m^{\ast}\in M^{\ast}$, $h\in H, h^{\ast}\in H^{\ast}$.
The cross symmetry is given by 
\be
\ps M^{\ast}, N, ((h^{\ast} \ot m^{\ast}) \ot (n \ot k)) := 
\Sigma \ k\1n \ot \ps \cu^{\ast}, \cv, ^C (m^{\ast} \ot n) \ot 
h^{\ast}\2n,
\ee
where $M^{\ast} := H^{\ast} \ot \cu^{\ast}$, $N := \cv \ot H$,
$\cu^{\ast} , \cv$ are objects in $\cp$, $m^{\ast}\in M^{\ast}$,
$n\in N$, $h^{\ast}\in H^{\ast}, k\in H$ and 
$\tg (h^{\ast}\ot k) := \Sigma \ k\1n \ot h^{\ast}\2n$.
\section{Hermitian Wick algebras and Fock space representation}
Let $\ca$ and $\cb$ be two unital and associative algebras over a 
field
$k$. The multiplication in these algebras is denoted by $\m \ca,$ and 
$\m \cb,$, respectively. Note that all tensor products are taken over
tha basic field $k$. Let us recall briefly the concept of a twisted 
product of algebras \cite{csv,bma}.
A $k$-linear mapping $\tau : \cb \ot \ca \ra \ca \ot \cb$ such that
$\tau (b\ot 1) = 1 \ot b, \:\:\: \tau (1 \ot a) = a \ot 1$ and
\be
\ba{l}
\tau \circ (id_{\cb} \ot m_{\ca}) = (m_{\ca} \ot id_{\cb}) 
\circ (id_{\ca} \ot \tau) \circ (\tau \ot id_{\ca})\\
\tau \circ (m_{\cb} \ot id_{\ca}) = (id_{\ca} \ot m_{\cb}) 
\circ (\tau \ot id_{\cb}) \circ (id_{\cb} \ot \tau)
\ea
\ee
is called a {\it twist}. We use the following standard notation
for the twist $\tau$
\be
\tau (b \ot a) := \Sigma \ a_{(i)} \ot b_{(i)}
\ee
It is known \cite{csv} that the multiplication $m_{\tau}$ defined 
on the tensor product $\ca \ot \cb$ by the formula
\be
\ba{c}
m_{\tau} := (m_A \ot m_B) \circ (id_A \ot \tau \ot id_B)
\label{mul}
\ea
\ee
is associative if and only if $\tau$ is a twist.
Let $\tau : \cb\ot\ca\rightarrow\ca\ot\cb$
be a twist. The tensor product of algebras $\ca\ot\cb$ 
equipped with the associative multiplication $m_{\tau}$ is said
to be a {\it twisted product} of algebras $\ca$ and $\cb$ with 
respect to the twist $\tau$ and is denoted by $\ca \tw_{\tau} \cb$.
Let $\ca$ be an associative algebra.
An algebra $\cb$ is said to be {\it conjugated} to the algebra
$\ca$ if there is an antilinear isomorphism 
$(-)^{\ast} : \ca \ra \cb$ such that
\be
(ab)^{\ast} = b^{\ast}a^{\ast},
\ee
where $a, b \in \ca$ and $a^{\ast}, b^{\ast}$ are their images 
under the isomorphism $(-)^{\ast}$. 
The inverse isomorphism $\ca^{\ast} \ra \ca$ will be denoted 
by the same symbol, i.e.
\be
(a^{\ast})^{\ast} = a.
\ee
If $\ca$ is an algebra, then the conjugated algebra will be denoted
by $\ca^{\ast}$. 
A twist $\tau : \ca^{\ast} \ot \ca \ra \ca \ot \ca^{\ast}$
such that
\be
(\tau (b^{\ast} \ot a))^{\ast} = \tau (a^{\ast} \ot b)
\ee
is called a $\ast$-twist.
Let $\tau : \ca^{\ast} \ot \ca \ra \ca \ot \ca^{\ast}$ be 
a $\ast$--twist. We define new algebra $\cw_{\tau}(\ca)$ as 
a twisted product of $\ca$ and $\ca^{\ast}$ with respect to 
the twist $\tau$, i.e.
\be
\ba{c}
\cw_{\tau}(\ca) := \ca \tw_{\tau} \ca^{\ast}.
\label{wck}
\ea
\ee
The $\ast$-operation in $\cw_{\tau}(\ca)$ is given by the relation
\be
\ba{c}
(a \ot b^{\ast})^{\ast} := (b \ot a^{\ast})
\label{hwa}
\ea
\ee
for $a, b \in \ca$. The algebra $\cw_{\tau}(\ca)$ defined by the 
relation (\ref{wck}) and equipped with $\ast$-operation 
(\ref{hwa}) is said to be a {\it Hermitian Wick} algebra. Let 
$\cp^{ext} := _{H^{\ast}} ^{H^{\ast}}\cp _H^H$ 
ba the category of right $H$--Hopf modules and left 
$H^{\ast}$-Hopf modules. 
Then there is a pair of conjugated algebras 
$\ca$ and $\ca^{\ast}$ in 
the category $\cp^{ext}$, the cross symmetry 
$\ps \ca^{\ast}, \ca, : \ca^{\ast}\ot\ca\ra\ca\ot\ca^{\ast}$,
the twist 
\be
\tau^{\Psi}_{g} := \ps \ca^{\ast}, \ca, + g^{ext},
\ee
and the corresponding Wick algebra 
$\cw_{\tau}(\ca) := \ca \tw_{\tau^{\Psi}_g} \ca^{\ast}$. 
Let $\ch$ be a $k$-vector space. We denote by $L(\ch)$ the algebra
of linear operators acting on $\ch$.\\
\bem
{\bf Theorem:} Let $\cw \equiv \ca \tt \tau, \ca^{\ast}$  
be a Hermitian Wick algebra. If $\pi_{\ca} : \ca \ra L(\ch)$ 
is a representation of the algebra $\ca$, such that we have 
the relation 
\be
\ba{c}
(\pi_{\ca}(b))^{\ast} \pi_{\ca}(a) = 
\pi_{\ca}(a_{(1)}) (\pi_{\ca}(b_{(2)}))^{\ast},
\label{wre}
\ea
\ee
then there is a representation $\pi_{\cw} : \cw \ra L(\ch)$ of the 
algebra $\cw$.\\
\eem
{\bf Proof:} If $\pi_{\ca} : \ca \ra L(\ch)$ is a representation 
of the algebra $\ca$, then the representation 
$\pi_{\ca^{\ast}} : \ca^{\ast} \ra L(\ch)$ of the conjugated 
algebra $\ca^{\ast}$ is given by the formula
\be
\pi_{\ca^{\ast}} (a^{\ast}) := (\pi_{\ca}(a))^{\ast}
\ee
for all $a \in \ca$. The existence of the representation 
$\pi_{\cw}$ of the algebra $\cw$ follows immediately from
the condition (\ref{wre}). 
\hfill $\Box$\\
For creation and annihilation operators we use the notation
\be
\pi_{\ca} (x^i) := a_{x^i}, \;\;\; 
\pi_{\ca^{\ast}} (x^{\ast i}) := a_{x^{\ast i}},
\ee
where $x^i$ and $x^{\ast i}$ are generators of $\ca$ and 
$\ca^{\ast}$, respectively, $i=1,\ldots ,N$.
We define creation operators for our model as multiplication
in the algebra $\ca$
\be
a_{x^i} s := m(x^i \ot s), \quad \mbox{for} \quad x^i, s \in \ca.
\ee
For the ground state and annihilation operators we assume that
\be
\langle 0|0 \rangle = 0, \quad a_{s^{\ast}} |0\rangle = 0 
\quad \mbox{for} \quad s^{\ast} \in \ca^{\ast},
\ee
and
\be
a_{x^{\ast i}} x^j := g^{ext} (x^{\ast i} \ot x^j).
\ee
It follows immediately from the relation (\ref{wre}) that
we have the following commutation relations
\be 
a_{x^{\ast i}} a_{x^i} - \Sigma \ 
a_{x^i \1n} a_{x^{\ast i}\2n} 
= g^{ext},
\ee
where $\ps \ca^{\ast}, \ca, (x^{\ast i}\ot x^j) 
= \Sigma \ x^j \1n \ot x^{\ast i}\2n$. 

Note that there is a probem with the existence of some 
additional relations \cite{RM,m10,ral}. For the desription 
of physical model we need the well--defined scalar product 
on the algebra of states $\ca$. For the the problem of the
existence of such scalar product and examples see the papers 
\cite{sin,wmq} and also \cite{qsym,qstat}. 
By the construction of the extended category we can classify
all possible interactions of particles with an external field.

\section*{Acknowledgments}  
The work is partially sponsored by Polish Committee 
for Scientific Research (KBN) under Grant 2P03B130.12.
The author would like to thank the organizer of SSPCM98
for his kind invitation to Zajaczkowo. The author would
like also to thank his own wife, A. Borowiec and all other
persons for the help which mades his participation in
the symposium possible.

\end{document}